\documentclass[reqno,11pt]{amsart}
\usepackage{amssymb}

\usepackage[english]{babel}
\usepackage[latin1]{inputenc}
\usepackage[T1]{fontenc}
\usepackage{latexsym,amsfonts}
\usepackage{graphics}
\usepackage{hhline}
\usepackage{mathrsfs}
\usepackage{amsmath}
\usepackage{amsthm}
\usepackage{pifont}
\usepackage{varioref}
\usepackage{textcomp}
\usepackage{lmodern}
\usepackage{euscript}
\usepackage[bookmarks=false]{hyperref}

\makeatletter\@addtoreset{equation}{section}

\newtheorem {theorem}{Theorem}[section]
\newtheorem {definition}[theorem]{Definition}
\newtheorem {remark}[theorem]{Remark}
\newtheorem {proposition}[theorem]{Proposition}

\newtheorem {lemma}[theorem]{Lemma}

\newcommand{\fin}{\hfill $\square$}

\begin{document}
\title[]{Coherent states quantization for generalized Bargmann spaces with formulae
for their attached Berezin transforms in terms of the Laplacian on $%
\mathbb{C}^{n}$}
\author{Zouhaïr MOUAYN}

\address{Department of Mathematics, Faculty of Sciences and Technics (M'Ghila) \\
Sultan Moulay Slimane University, PO. Box 523, B\'{e}ni Mellal, Morocco}
\email{mouayn@fstb.ac.ma}
\maketitle

\begin{abstract}
While dealing with a class of generalized Bargmann spaces, we rederive their
reproducing kernels from the knowledge of an orthonormal basis by using an
addition formula for Laguerre polynomials involving the disk polynomials. We
construct for each of \ these spaces a set of coherent states to apply a
coherent states quantization method. This provides us with another way to
recover the Berezin transforms attached to these spaces. Finally, two new
formulae representing these transforms as functions of the Euclidean
Laplacian are established and a possible physics direction for the
application of such formulae is discussed.
\end{abstract}

\section{Introduction}

The Berezin transform introduced in \cite{1} for certain classical symmetric
domains in $\mathbb{C}^{n}$ is a transform linking the Berezin symbols and
the symbols for Toeplitz operators. It is present in the study of the
correspondence principle. The formula representing the Berezin transform as
a function of the Laplace-Beltrami operator plays a key role in the Berezin
quantization \cite{2}.

In this paper we are concerned with the Berezin transforms associated with
\begin{equation}
A_{m}^{2}(\mathbb{C}^{n})=\left\{ \psi \in L^{2}(\mathbb{C}%
^{n};e^{-|z|^{2}}d\mu );\qquad \widetilde{\Delta }\psi =m\psi \right\},
\label{1.1}
\end{equation}
the eigenspaces (\cite{3}) of the second order differential operator
\begin{equation}
\widetilde{\Delta }=-\sum\limits_{j=1}^{n}\frac{\partial ^{2}}{\partial
z_{j}\partial \overline{z}_{j}}+\sum\limits_{j=1}^{n}\overline{z}_{j}\frac{%
\partial }{\partial \overline{z}_{j}} \label{1.2}
\end{equation}
corresponding to the eigenvalues $m=0,1,2,\ldots $. \ Here, $d\mu $ denotes
the Lebesgue measure on $\mathbb{C}^{n}$. The eigenspaces in \eqref{1.1} are
called generalized Bargmann spaces and are reproducing kernel Hilbert spaces
with reproducing kernels given by (\cite{3}):
\begin{equation}
K_{m}\left( z,w\right) :=\pi ^{-n}e^{\left\langle z,w\right\rangle
}L_{m}^{\left( n-1\right) }\left( \left| z-w\right| ^{2}\right) ,\qquad
z,w\in \mathbb{C}^{n},  \label{1.3}
\end{equation}
where $L_{k}^{(\alpha )}(\cdot )$ is the Laguerre polynomial \cite{4}. The
associated Berezin transform was obtained via the well known formalism of
Toeplitz operators by the convolution product over $\mathbb{C}^{n}$ as (\cite
{5,6}):
\begin{equation}
B_{m}[\phi ](z):=\left( h_{m}*\phi \right) (z), \qquad \phi \in L^{2}(\mathbb{C}^{n},d\mu ),  \label{1.4}
\end{equation}
where
\begin{equation}
h_{m}\left( w\right) :=\frac{\pi ^{-n}m!}{(n)_{m}}e^{-|w|^{2}}\left(
L_{m}^{(n-1)}\left( |w|^{2}\right) \right) ^{2},\qquad w\in \mathbb{C}^{n}.
\label{1.5}
\end{equation}

Here, our aim is first to present a direct proof for obtaining the
reproducing kernel in Eq. \eqref{1.3} starting from the knowledge of an
orthonormal basis of the eigenspace in \eqref{1.1} and by exploiting an
addition formula for Laguerre polynomials involving the disk polynomials
\cite{17}. Next, we construct for each of \ the eigenspaces in \eqref{1.1} a
set of coherent states by following a generalized formalism \cite{8}. in
order to apply a coherent states quantization method. This provides us with
another way to recover the Berezin transforms in \eqref{1.4} attached to the
generalized Bargmann spaces under consideration. Finally, direct
calculations of the Fourier transform of the function in \eqref{1.5} enables
us to present two other formulae expressing the transforms in \eqref{1.4} as
functions of the Euclidean Laplacian. These two formulae together with a
first on \cite{6} could be of help in physics problems having a close
analogy with the diamagnetism of spinless bose systems \cite{9}.

This paper is organized as follows. In Section 2 we recall briefly the
formalism of coherent states quantization we will be using. Section 3 deals
with some needed facts on the generalized Bargmann spaces. In Section 4, we
construct for each of these spaces a set of coherent states and we apply the
corresponding quantization scheme in order to recover the Berezin
transforms. In Section 5, we present two other formulae expressing these
Berezin transforms as functions of the Laplacian in the Euclidean complex $n$-space. Section 6 is devoted to some concluding remarks.

\section{Coherent states quantization}

Coherent states are mathematical tools which provide a close connection
between classical and quantum formalisms. In general, they are a specific
overcomplete set of vectors in a Hilbert space satisfying a certain
resolution of the identity condition. Here, we adopt the prototypical model
of coherent states presented in \cite{8} and described as follows. Let $%
X=\left\{ x\mid x\in X\right\} $ be a set equipped with a measure $d\nu $
and $L^{2}(X,d\nu )$ the Hilbert space of $d\nu -$\ square integrable
functions $f(x)$ on $X$. Let $\mathcal{A}^{2}\subset L^{2}(X,d\nu )$ be a
subspace with an orthonormal basis $\{\Phi _{k}\}_{k=0}^{\infty }$ such that
\begin{equation}
\mathcal{N}(x):=\sum_{k=0}^{\infty }\left| \Phi _{k}(x)\right| ^{2}<+\infty
,\quad x\in X.  \label{2.1}
\end{equation}
Let $\mathcal{H}$ be another (functional) Hilbert space with $\dim \mathcal{H%
}=\dim \mathcal{A}^{2}$ and $\left\{ \varphi _{k}\right\} _{k=0}^{\infty }$
is a given orthonormal basis of $H$. Then consider the family of states $%
\left\{ |x>\right\} _{x\in X}$ in $\mathcal{H}$, through the following
linear superpositions:
\begin{equation}
|x>:=\left( N(x)\right) ^{-\frac{1}{2}}\sum_{k=0}^{\infty }\Phi
_{k}(x)\varphi _{k}.  \label{2.2}
\end{equation}
These \textit{coherent states} obey the normalization condition
\begin{equation}
<x\mid x>_{\mathcal{H}}=1  \label{2.3}
\end{equation}
and the following resolution of the unity in $\mathcal{H}$
\begin{equation}
\mathbf{1}_{\mathcal{H}}=\int_{X}|x><x|\mathcal{N}(x)d\mu (x),  \label{2.4}
\end{equation}
which is expressed in terms of Dirac's bra-ket notation $|x><x|$ meaning the
rank-one-operator $\varphi \mapsto <\varphi \mid x>_{\mathcal{H}}.|x>$, $
\varphi \in \mathcal{H}$.

The choice of the Hilbert space $\mathcal{H}$ defines in fact a quantization
of the space $X$ by the coherent states in \eqref{2.2}, via the inclusion
map $X \ni x\mapsto |x>\in \mathcal{H}$ and the property \eqref{2.4} is
crucial in setting the bridge between the classical and the quantum worlds.
It encodes the quality of being canonical quantizers along a guideline
established by Klauder \cite{10} and Berezin \cite{1}. This Klauder-Berezin
coherent state quantization, also named anti-Wick quantization or Toeplitz
quantization \cite{11} by many authors, consists in associating to a
classical observable that is a function $f(x)$ on $X$ having specific
properties the operator-valued integral
\begin{equation}
A_{f}:=\int_{X}|x><x|f(x)\mathcal{N}(x)d\mu (x).  \label{2.5}
\end{equation}
The function $f(x)\equiv $ $\widehat{A}_{f}(x)$ is called upper (or
contravariant) symbol of the operator $A_{f}$ and is nonunique in general.
On the other hand, the expectation value $<x\mid $ $A_{f}\mid x>$ of $A_{f}$
\ with respect to the set of coherent states $\left\{ |x>\right\} _{x\in X}$
is called lower (or covariant) symbol of $A_{f}.$ Finally, associating to
the classical observable $f(x)$ the obtained mean value $<x\mid $ $A_{f}\mid
x>,$ we get the Berezin transform of this observable. That is,
\begin{equation}
B\left[ f\right] (x):=<x\mid A_{f}\mid x>,\qquad x\in X.  \label{2.6}
\end{equation}
For all aspect of the theory of coherent states and their genesis, we refer
to the survey \cite{12} by V.V. Dodonov or the recent book \cite{8} of J.P.
Gazeau.

\section{\noindent The generalized Bargmann spaces $A_{m}^{2}(\mathbb{C}%
^{n}) $}

In \cite{3}, we have introduced \ a class of generalized Bargmann spaces
indexed by integer numbers $m=0,1,2,\cdots ,$ as eigenspaces of a single
elliptic second-order differential operator as
\begin{equation}
A_{m}^{2}(\mathbb{C}^{n}):=\left\{ \psi \in L^{2}(\mathbb{C}%
^{n};e^{-|z|^{2}}d\mu );\qquad \widetilde{\Delta }\psi =m\psi \right\} ,
\label{3.1}
\end{equation}
where
\begin{equation}
\widetilde{\Delta }=-\sum\limits_{j=1}^{n}\frac{\partial ^{2}}{\partial
z_{j}\partial \overline{z}_{j}}+\sum\limits_{j=1}^{n}\overline{z}_{j}\frac{%
\partial }{\partial \overline{z}_{j}}.  \label{3.2}
\end{equation}
The operator $\widetilde{\Delta }$ is considered with $C_{0}^{\infty }(%
\mathbb{C}^{n})$ as its regular domain in the Hilbert space $L^{2}(\mathbb{C}%
^{n};e^{-|z|^{2}}d\mu )$ of $e^{-|z|^{2}}d\mu $-square integrable
complex-valued functions on $\mathbb{C}^{n}.$ Its concrete $L^{2}$ spectral
theory have been discussed in \cite{13}. Actually, for $m=0,$ the eigenspace
$A_{0}^{2}(\mathbb{C}^{n})$ turns out to be the realization by harmonic
functions with respect to $\widetilde{\Delta }$ of the classical Bargmann
space whose elements are the entire functions in $L^{2}(\mathbb{C}%
^{n};e^{-|z|^{2}}d\mu )$ see \cite{3}.

Now, for general $m=0,1,2,\cdots $, a complete description of the expansion
of elements $f\in A_{m}^{2}(\mathbb{C}^{n})$ in terms of the appropriate
Fourier series in $\mathbb{C}^{n}$ have been given \cite{3}. Precisely, a
function $f:\mathbb{C}^{n}\rightarrow \mathbb{C}$ belongs to $A_{m}^{2}(%
\mathbb{C}^{n})$ if and only if it can be expanded in the form
\begin{equation}
f(z) =\sum\limits_{p=0}^{+\infty }\sum\limits_{q=0}^{m} \quad
_{1}F_{1}\left( -m+q,n+p+q,\rho ^{2}\right) \rho ^{p+q}\gamma _{p,q}\cdot
h_{p,q}\left( \omega \right)  \label{3.3}
\end{equation}
in $C^{\infty }\left( \mathbb{C}^{n}\right) ,$ $z=\rho \omega ,\rho >0,$ $%
\left| \omega \right| =1$, $_{1}F_{1}$ is the confluent hypergeometric
function \cite{4}, $\gamma _{p,q}=\left( \gamma _{p,q,j}\right) \in %
\mathbb{C}^{d\left( n;p,q\right) }$ are coefficients such that
\begin{equation}
\sum\limits_{p=0}^{+\infty }\sum\limits_{q=0}^{m}\left( m-q\right) !\left(
p+q+n-1\right) !\Gamma \left( n+p+q\right) \frac{\left| \gamma _{p,q}\right|
^{2}}{2\Gamma \left( n+m+p\right) }<+\infty  \label{3.4}
\end{equation}
and
\begin{equation}
d\left( n;p,q\right) :=\frac{\left( p+q+n-1\right) \left( p+n-2\right)
!\left( q+n-2\right) !}{p!q!\left( n-1\right) !\left( n-2\right) !}
\label{3.5}
\end{equation}
is the dimension of the space $H\left( p,q\right) $ of restrictions to the
unit sphere $\mathbb{S}^{2n-1}=\left\{ \omega \in \mathbb{C}^{n},\left|
\omega \right| =1\right\} $ of harmonic polynomials $h(z) $ on $\mathbb{C}%
^{n}$, which are homogeneous of degree $p$ in $z$ and degree $q$ in $%
\overline{z}$ (see \cite{14}). The notation ``$\cdot "$ in \eqref{3.3} means
the following finite sum
\begin{equation}
\gamma _{p,q}\cdot h_{p,q}\left( \omega \right) :=\sum\limits_{j=0}^{d\left(
n,p,q\right) }\gamma _{p,q,j}.h_{p,q}^{j}\left( \omega \right) ,  \label{3.6}
\end{equation}
where $\left\{ h_{p,q}^{j}\right\} _{1\leq j\leq d\left( n;p,q\right) }$ is
an orthonormal basis of $H\left( p,q\right) .$

Now, from \eqref{3.3} and using the relation (\cite{3}):
\begin{equation}
_{1}F_{1}\left( -k,\alpha ;u\right) =\frac{k!\Gamma \left( \alpha +1\right)
}{\Gamma \left( \alpha +k+1\right) }L_{k}^{(\alpha )}(u)  \label{3.7}
\end{equation}
for $k$ being a positive integer, an orthonormal basis in the space $%
A_{m}^{2}(\mathbb{C}^{n})$ can be written explicitly in terms of the
Laguerre polynomials $L_{k}^{(\alpha )}(.)$ and the polynomials $%
h_{p,q}^{j}(z,\overline{z})$ as \cite{6}:
\begin{equation}
\Phi _{j,p,q}^{m}(z):=\left( \frac{2\left( m-q\right) !}{\Gamma \left(
n+m+p\right) }\right) ^{\frac{1}{2}}L_{m-q}^{\left( n+p+q-1\right) }\left(
|z|^{2}\right) h_{p,q}^{j}\left( z,\overline{z}\right)   \label{3.8}
\end{equation}
for varying $p=0,1,2,\ldots ;q=0,1,\ldots ,m$ and $j=1,\ldots ,d\left(
n;p,q\right) $. These functions possess nice properties and represent a
principal tool in the present work. For instance, we have to  check by hand
the following fact.

\begin{lemma}
\label{Lemma3.1} The functions in \eqref{3.8} satisfy
\begin{equation}
\sum\limits_{p=0}^{+\infty }\sum\limits_{q=0}^{m}\sum_{j=1}^{d\left(
n,p,q\right) }\Phi _{j,p,q}^{m}(z)\overline{\Phi _{j,p,q}^{m}(w)}=\pi
^{-n}e^{\left\langle z,w\right\rangle }L_{m}^{\left( n-1\right) }\left(
\left| z-w\right| ^{2}\right)   \label{3.9}
\end{equation}
for all $z,w\in \mathbb{C}^{n}.$
\end{lemma}

\noindent \textit{Proof.} Denoting the sum in the left hand side of %
\eqref{3.9} by $S_{n,m}^{z,w}$ and replacing the functions $\Phi
_{j,p,q}^{m}(z)$ by their expressions in \eqref{3.8}, we obtain that
\begin{align}
S_{n,m}^{z,w}& :=\sum\limits_{\substack{ 1\leq j\leq d\left( n,p,q\right) \\
0\leq q\leq m,0\leq p<+\infty }}\Phi _{j,p,q}^{m}(z)\overline{\Phi
_{j,p,q}^{m}(w)}  \label{3.10} \\
& =\sum\limits_{\substack{ 1\leq j\leq d\left( n,p,q\right) \\ 0\leq q\leq
m,0\leq p<+\infty }}\frac{2\left( m-q\right) !}{\Gamma \left( n+m+p\right) }
L_{m-q}^{\left( n+p+q-1\right) }\left( |z|^{2}\right) L_{m-q}^{\left(
n+p+q-1\right) }\left( |w|^{2}\right) h_{p,q}^{j}\left( z,\overline{z}%
\right) \overline{h_{p,q}^{j}\left( w,\overline{w}\right) }  \label{3.11} \\
& =\sum\limits_{\substack{ 0\leq q\leq m, \\ 0\leq p<+\infty }}\frac{2\left(
m-q\right) !}{\Gamma \left( n+m+p\right) }L_{m-q}^{\left( n+p+q-1\right)
}\left( |z|^{2}\right) L_{m-q}^{\left( n+p+q-1\right) }\left( |w|^{2}\right)
\frak{G}_{n;p,q}^{z,w},  \label{3.12}
\end{align}
where
\begin{align}
\frak{G}_{n;p,q}^{z,w}& :=\sum\limits_{1\leq j\leq d\left( n,p,q\right)
}h_{p,q}^{j}\left( z,\overline{z}\right) \overline{h_{p,q}^{j}\left( w,%
\overline{w}\right) }  \label{3.13} \\
& =\left( \left| z\right| \left| w\right| \right) ^{p+q}\sum\limits_{1\leq
j\leq d\left( n,p,q\right) }h_{p,q}^{j}\left( \frac{z}{\left| z\right| },%
\frac{\overline{z}}{\left| z\right| }\right) \overline{h_{p,q}^{j}\left(
\frac{w}{\left| w\right| },\frac{\overline{w}}{\left| w\right| }\right) }.
\label{3.14}
\end{align}
To calculate the finite sum in \eqref{3.13}, we make use of the
formula (\cite{15}):
\begin{align}
\sum\limits_{1\leq j\leq d\left( n,p,q\right) }& h_{p,q}^{j}\left( \frac{z}{%
\left| z\right| },\frac{\overline{z}}{\left| z\right| }\right) \overline{%
h_{p,q}^{j}\left( \frac{w}{\left| w\right| },\frac{\overline{w}}{\left|
w\right| }\right) }=\frac{\Gamma \left( n\right) }{2\pi ^{n}}d\left(
n,p,q\right) \left( P_{\min \left( p,q\right) }^{\left( n-2,\left|
p-q\right| \right) }\left( 1\right) \right) ^{-1}  \label{3.15} \\
& \times \left| \left\langle \frac{z}{\left| z\right| },\frac{w}{\left|
w\right| }\right\rangle \right| ^{\left| p-q\right| }e^{i\left( p-q\right)
\arg \left\langle \frac{z}{\left| z\right| },\frac{w}{\left| w\right| }%
\right\rangle }.P_{\min \left( p,q\right) }^{\left( n-2,\left| p-q\right|
\right) }\left( 2\left| \left\langle \frac{z}{\left| z\right| },\frac{w}{%
\left| w\right| }\right\rangle \right| ^{2}-1\right) ,  \notag
\end{align}
where $P_{k}^{\left( \alpha ,\beta \right) }(.)$ is the Jacobi polynomial
\cite{4}. Now, making appeal to the normalized Jacobi polynomial
\begin{equation}
R_{k}^{\left( \alpha ,\beta \right) }\left( u\right) :=\frac{P_{k}^{\left( \alpha
,\beta \right) }(.)}{P_{k}^{\left( \alpha ,\beta \right) }\left( 1\right)}
\label{3.16}
\end{equation}
together with the disk polynomials that were first studied by Zernike and
Brinkman \cite{16} and are given by
\begin{equation}
R_{p,q}^{\gamma }\left( \xi \right) :=\left| \xi \right| ^{\left| p-q\right|
}e^{i\left( p-q\right) \arg \xi }R_{\min \left( p,q\right) }^{\left( \gamma
,\left| p-q\right| \right) }\left( 2\left| \xi \right| ^{2}-1\right) ,
\label{3.17}
\end{equation}
where $R_{\min \left( p,q\right) }^{\left( \gamma ,\left| p-q\right| \right)
}(.)$ is defined according to \eqref{3.16} for the parameter $\gamma :=n-2,$
the finite sum takes the form
\begin{equation}
\frak{G}_{n;p,q}^{z,w}=\left( 2\pi ^{n}\right) ^{-1}\Gamma \left( n\right)
d\left( n,p,q\right) \left( \left| z\right| \left| w\right| \right)
^{p+q}R_{p,q}^{n-2}\left( \left\langle \frac{z}{\left| z\right| },\frac{w}{%
\left| w\right| }\right\rangle \right) .  \label{3.18}
\end{equation}
Returning back to \eqref{3.12} and inserting \eqref{3.18}, then the sum $%
S_{n,m}^{z,w}$ takes the form
\begin{align}
S_{n,m}^{z,w}& =\sum\limits_{\substack{ 0\leq q\leq m, \\ 0\leq p<+\infty }}%
\frac{\left( m-q\right) !\Gamma \left( n\right) d\left( n,p,q\right) }{%
\Gamma \left( n+m+p\right) \pi ^{n}}\left( \left| z\right| \left| w\right|
\right) ^{p+q}  \label{3.19} \\
& \times L_{m-q}^{\left( n+p+q-1\right) }\left( |z|^{2}\right)
L_{m-q}^{\left( n+p+q-1\right) }\left( |w|^{2}\right) R_{p,q}^{n-2}\left(
\left\langle \frac{z}{\left| z\right| },\frac{w}{\left| w\right| }%
\right\rangle \right) .  \notag
\end{align}
Now, we use the slightly different function $\mathcal{L}_{k}^{\left( \alpha
\right) }\left( u\right) $ for the Laguerre polynomial, which is such that
\begin{equation}
L_{k}^{\left( \alpha \right) }\left( u\right) =\frac{\Gamma \left( k+\alpha
+1\right) }{k!\Gamma \left( \alpha +1\right) }e^{\frac{1}{2}u}\mathcal{L}%
_{k}^{\left( \alpha \right) }\left( u\right)   \label{3.20}
\end{equation}
for $\alpha =n-1+p+q$, $u=|z|^{2}$, $\ u=|w|^{2}$ and $k=m-q$ to rewrite %
\eqref{3.19} as
\begin{align}
S_{n,m}^{z,w}& =\pi ^{-n}e^{\frac{1}{2}|z|^{2}+\frac{1}{2}\left| w\right|
^{2}}\sum\limits_{\substack{ 0\leq q\leq m, \\ 0\leq p<+\infty }}\frac{%
\Gamma \left( n\right) d\left( n,p,q\right) \Gamma \left( m+n+p\right) }{%
\left( m-q\right) !\Gamma ^{2}\left( n+p+q\right) }  \label{3.21} \\
& \times \left( \left| z\right| \left| w\right| \right) ^{p+q}\mathcal{L}%
_{m-q}^{\left( n-1+p+q\right) }\left( |z|^{2}\right) \mathcal{L}%
_{m-q}^{\left( n-1+p+q\right) }\left( |w|^{2\text{ }}\right)
R_{p,q}^{n-2}\left( \left\langle \frac{z}{\left| z\right| },\frac{w}{\left|
w\right| }\right\rangle \right) .  \notag
\end{align}
Next, we use the explicit expression of the dimension $d\left( n,p,q\right) $
in \eqref{3.5} to write the coefficient
\begin{equation}
C_{n,p,q}:=\frac{\Gamma \left( n\right) d\left( n,p,q\right) \Gamma \left(
m+n+p\right) }{\left( m-q\right) !\Gamma ^{2}\left( n+p+q\right) }
\label{3.22}
\end{equation}
in the following form
\begin{equation}
C_{n,p,q}=\frac{1}{\left( n-2\right) !}\left[ \frac{\Gamma \left( n+m\right)
}{m!\left( n-1\right) }\right] \frac{\sigma }{\sigma +p+q}\left(
\begin{array}{c}
m \\
q
\end{array}
\right) \frac{\left( \sigma +m+1\right) _{p}}{p!\left( \sigma +q\right)
_{p}\left( \sigma +p\right) _{q}},  \label{3.23}
\end{equation}
where $\sigma =n-1$ and $\left( a\right) _{k}$ denote the Pochammer symbol.
Therefore, the sum in \eqref{3.21} can be presented as
\begin{align}
S_{n,m}^{z,w}& =\frac{\Gamma \left( n+m\right) e^{\frac{1}{2}(|z|^{2}+\left|
w\right| ^{2})}}{\pi ^{n}m!\left( n-1\right) !}\sum\limits_{p=0}^{+\infty
}\sum\limits_{q=0}^{m}\frac{\sigma }{\sigma +p+q}\left(
\begin{array}{c}
m \\
q
\end{array}
\right) \frac{\left( \sigma +m+1\right) _{p}}{p!\left( \sigma +q\right)
_{p}\left( \sigma +p\right) _{q}}  \label{3.24} \\
& \left| z\right| ^{p+q}\mathcal{L}_{m-q}^{\left( \sigma +p+q\right) }\left(
|z|^{2}\right) \left| w\right| ^{p+q}\mathcal{L}_{m-q}^{\left( \sigma
+p+q\right) }\left( |w|^{2\text{ }}\right) R_{p,q}^{\sigma -1}\left( \rho
e^{i\theta }\right) ,  \notag
\end{align}
where $\left\langle \frac{z}{\left| z\right| },\frac{w}{\left| w\right| }%
\right\rangle =\rho e^{i\theta }$.

We are now in position to apply the addition formula for Laguerre
polynomials due to Koornwinder \cite{17}:
\begin{align}
\exp \left( ixyr\sin \psi \right) & \mathcal{L}_{s}^{\left( \sigma \right)
}\left( x^{2}+y^{2}-2xy\cos \psi \right)  \label{3.25} \\
& =\sum\limits_{k=0}^{+\infty }\sum\limits_{l=0}^{s}\frac{\sigma }{\sigma
+k+l}\left(
\begin{array}{c}
s \\
l
\end{array}
\right) \frac{\left( \sigma +s+1\right) _{k}}{k!\left( \sigma +l\right)
_{k}\left( \sigma +k\right) _{l}}  \notag \\
& \times x^{k+l}\mathcal{L}_{s-l}^{\left( \sigma +k+l\right) }\left(
x^{2}\right) y^{k+l}\mathcal{L}_{s-l}^{\left( \sigma +k+l\right) }\left(
y^{2}\right) R_{k,l}^{\sigma -1}\left( re^{i\psi }\right) ,  \notag
\end{align}
where $x\geq 0,y\geq 0,0\leq r\leq 1,0\leq \psi <2\pi ,\sigma >0, \quad
s=0,1,2,\cdots , $ for $s=m,\sigma =n-1,re^{i\psi }=\left\langle \frac{z}{%
\left| z\right| },\frac{w}{\left| w\right| }\right\rangle ,x=\left| z\right|
$ and $y=\left| w\right| $. After computations, we arrive at
\begin{equation}
S_{n,m}^{z,w}=\pi ^{-n}e^{\left\langle z,w\right\rangle }L_{m}^{\left(
n-1\right) }\left( \left| z-w\right| ^{2}\right) .  \label{3.26}
\end{equation}
This ends the proof. \hfill $\square$\newline

By a general fact on reproducing kernels \cite{18}, Lemma \ref{Lemma3.1}
says that the knowledge of the explicit orthonormal basis in \eqref{3.8}
leads directly to the expression of the reproducing kernel in \eqref{1.2}
via calculations using properties of the spherical harmonics together with
the addition formula \eqref{3.25}.

\begin{remark}
The motion of charged particle in a constant uniform magnetic field in $%
\mathbb{R}^{2n}$ is described (in suitable units and up to additive
constant) by the Schrödinger operator
\begin{equation}
H_{\beta }:=-\frac{1}{4}\sum\limits_{j=1}^{n}\left( \partial _{x_{j}}+\beta
y_{j}\right) ^{2}+\left( \partial _{y_{j}}-i\beta x_{j}\right) ^{2}-\frac{n}{%
2}  \label{3.27}
\end{equation}
acting on $L^{2}\left( {}\mathbb{R}^{2n},d\mu \right) $, $\beta >0$ is a
constant proportional to the magnetic field strength. We identify the
Euclidean space ${}\mathbb{R}^{2n}$ with $\mathbb{C}^{n}$ in the usual way.
The operator $H_{\beta }$ in \eqref{3.27} can be represented by the operator
\begin{equation}
\widetilde{H}_{\beta }=e^{\frac{1}{2}\beta |z|^{2}}H_{\beta }e^{-\frac{1}{2}%
\beta |z|^{2}}.  \label{3.28}
\end{equation}
According to equation \eqref{3.28}, an arbitrary state $\phi $ of $%
L^{2}\left( \mathbb{{}R}^{2n},d\mu \right) $ \hspace{0in}is represented by
the function $Q\left[ \phi \right] $ of $L^{2}\left( \mathbb{C}%
^{n},e^{-|z|^{2}}d\mu \right) $ defined by
\begin{equation}
Q\left[ \phi \right] (z):=e^{\frac{1}{2}\beta |z|^{2}}\phi (z),z\in %
\mathbb{C}^{n}.  \label{3.29}
\end{equation}
The unitary map $Q$ in \eqref{3.29} is called a ground state transformation.
For $\beta =1$, the explicit expression for the operator in equation %
\eqref{3.28} turns out to be given by the operator $\widetilde{\Delta }$
introduced in equation \eqref{1.2}. That is, $\widetilde{H}_{1}=$ $%
\widetilde{\Delta }$ with eigenvalues corresponding to well known \textit{%
Landau levels} of the operator in \eqref{3.27}.
\end{remark}

\section{A coherent states quantization for $A_{m}^{2}(\mathbb{C}^{n})$}

Now, to adapt the definition \eqref{2.1} of coherent states for the context
of the generalized Bargmann spaces $A_{m}^{2}(\mathbb{C}^{n}),$ we first
list the following notations.

\begin{itemize}
\item[$\bullet $]  $\left( X,d\nu \right) :=(\mathbb{C}^{n},e^{-|z|^{2}}d\mu
),$ $d\nu :=e^{-|z|^{2}}d\mu $.

\item[$\bullet $]  $x\equiv z\in \mathbb{C}^{n}$.

\item[$\bullet $]  $A^{2}:=A_{m}^{2}(\mathbb{C}^{n})\subset L^{2}\left( %
\mathbb{C}^{n},e^{-|z|^{2}}d\mu \right) $.

\item[$\bullet $]  $\left\{ \Phi _{k}(x)\right\} \equiv \left\{ \Psi
_{j,p,q}^{m}(z)\right\} $ is the orthonormal basis of $A_{m}^{2}(\mathbb{C}%
^{n})$ in \eqref{3.8}.

\item[$\bullet $]  $\mathcal{N}(x)\equiv \mathcal{N}(z)$ is a normalization
factor.

\item[$\bullet $]  $\left\{ \varphi _{k}\right\} \equiv \left\{ \varphi
_{j,p,q}\right\} $ is an orthonormal basis of another (functional) Hilbert
space $\mathcal{H}$ having the same dimension $\left( \infty \right) $ of $%
A_{m}^{2}(\mathbb{C}^{n})$.
\end{itemize}

\begin{definition}
\label{Defintion4.1} For each fixed integer $m=0,1,2,\cdots ,$ a class of
the generalized coherent states associated with the space $A_{m}^{2}(C^{n})$
is defined according to \eqref{2.1} by the form
\begin{equation}
\phi _{z,m}\equiv \mid z,m>:=\left( \mathcal{N}(z)\right) ^{-\frac{1}{2}%
}\sum\limits_{\substack{ 1\leq j\leq d\left( n,p,q\right) \\ 0\leq q\leq
m,0\leq p<+\infty }}\Phi _{j,p,q}^{m}(z)\varphi _{j,p,q},  \label{4.1}
\end{equation}
where $\mathcal{N}(z)$ is a factor such that $<\phi _{z,m},\phi _{z,m}>_{%
\mathcal{H}}=1.$
\end{definition}

\begin{proposition}
\label{Proposition4.2.} The normalization factor in \eqref{4.1} is given by
\begin{equation}
\mathcal{N}(z)=\frac{\pi ^{-n}\Gamma \left( n+m\right) }{\Gamma \left(
m+1\right) \Gamma \left( n\right) }e^{\left\langle z,z\right\rangle }
\label{4.2}
\end{equation}
for every $z\in \mathbb{C}^{n}.$
\end{proposition}

\noindent \textit{Proof}. To calculate this factor, we start by writing the
condition
\begin{equation}
<\phi _{z,m},\phi _{z,m}>_{\mathcal{H}}=1.  \label{4.3}
\end{equation}
Eq. \eqref{4.3} is equivalent to
\begin{equation}
\left( \mathcal{N}(z) \right) ^{-1}\sum\limits_{p=0}^{+\infty
}\sum\limits_{q=0}^{m}\sum_{j=1}^{d\left( n,p,q\right) }\Phi _{j,p,q}^{m}(z)%
\overline{\Phi _{j,p,q}^{m}(z)}=1.  \label{4.4}
\end{equation}
Making use of Lemma \ref{Lemma3.1} for the particular case $z=w$, we get
that
\begin{equation}
\mathcal{N}(z) =\pi ^{-n}e^{\left\langle z,z\right\rangle }L_{m}^{\left(
n-1\right) }\left( 0\right) .  \label{4.5}
\end{equation}
Next, by the fact that (\cite{19})
\begin{equation}
L_{m}^{\left( n-1\right) }\left( 0\right) =\frac{\Gamma \left( n+m\right) }{%
\Gamma \left( m+1\right) \Gamma \left( n\right) },  \label{4.6}
\end{equation}
we arrive at the announced result. \hfill $\square$ \newline

Now, the states $\phi _{z,m}$\ $\equiv \mid z,m>$\ satisfy the resolution of
the identity
\begin{equation}
1_{\mathcal{H}}=\int_{\mathbb{C}^{n}}\mid z,m><z,m\mid \mathcal{N}(z) d\nu
(z)  \label{4.7}
\end{equation}
and with the help of them we can achieve the coherent states quantization
scheme described in Section 2 to rederive the Berezin transform $B_{m}$ in %
\eqref{1.4}  which was defined by the Toeplitz operators formalism in
previous works. For this let us first associate to any arbitrary function $%
\varphi \in L^{2}\left( \mathbb{C}^{n},d\mu \right) $ the operator-valued
integral
\begin{equation}
A_{\varphi }:=\int_{\mathbb{C}^{n}}|z,m><m,z|\varphi (z) \mathcal{N}(z)d\nu
(z).  \label{4.8}
\end{equation}
The function $\varphi (z) $ is a upper symbol of the operator $A_{\varphi }$. On the other hand, we need to calculate the expectation value
\begin{equation}
\mathbb{E}_{\left\{ \mid z,m>\right\} }\left( A_{\varphi }\right) :=<z,m\mid
A_{\varphi }\mid z,m>  \label{4.9}
\end{equation}
of $A_{\varphi }$ with respect to the set of coherent states $\left\{ \mid
z,m>\right\} _{z\in \mathbb{C}}$ defined in \eqref{4.1}. This will
constitute a lower symbol of the operator $A_{\varphi }$.

\begin{proposition}
\label{Proposition4.3} Let $\varphi \in L^{2}\left( \mathbb{C}^{n},d\mu
\right) $. Then, the expectation value in \eqref{4.9} has the following
expression
\begin{equation}
\mathbb{E}_{\left\{ \mid z,m>\right\} }\left( A_{\varphi }\right) =\frac{m!}{%
(n)_{m}\pi ^{n}}\int_{\mathbb{C}^{n}}e^{-|z-w|^{2}}\left( L_{m}^{\left(
n-1\right) }\left( |z-w|^{2}\right) \right) ^{2}\varphi (w)d\mu (w)
\label{4.10}
\end{equation}
for every $z\in \mathbb{C}^{n}.$
\end{proposition}

\noindent \textit{Proof. } We first write the action of the operator $%
A_{\varphi }$ in \eqref{4.8} on an arbitrary coherent state $\mid z,m>$ in
terms of Dirac's bra-ket notation as
\begin{equation}
A_{\varphi }\mid z,m>=\int_{\mathbb{C}^{n}}|w,m><w,m|z,m>\varphi \left(
w\right) \mathcal{N}(w)e^{-|w|^{2}}d\mu \left( w\right) \label{4.11}
\end{equation}
Therefore, the expectation value reads
\begin{align}
<z,m\mid A_{\varphi }\mid z,m>& =\int_{\mathbb{C}^{n}}<z,m|w,m><w,m|z,m>%
\varphi \left( w\right) \mathcal{N}(w)e^{-|w|^{2}}d\mu \left( w\right)
\label{4.12} \\
& =\int_{\mathbb{C}^{n}}<z,m|w,m>\overline{<z,m|w,m>}\varphi \left( w\right)
\mathcal{N}(w)e^{-|w|^{2}}d\mu \left( w\right)   \label{4.13} \\
& =\int_{\mathbb{C}^{n}}\left| <z,m|w,m>\right| ^{2}\varphi \left( w\right)
\mathcal{N}(w)e^{-|w|^{2}}d\mu \left( w\right) .  \label{4.14}
\end{align}
Note that we have used the fact that $d\nu \left( w\right)
:=e^{-|w|^{2}}d\mu \left( w\right)$. Now, we need to evaluate the quantity $%
\left| <z,m|w,m>\right| ^{2}$ in \eqref{4.14}. For this, we write
the scalar product as
\begin{equation}
<z,m|w,m>=\sum\limits_{p=0}^{+\infty
}\sum\limits_{q=0}^{m}\sum_{j=1}^{d\left( n,p,q\right)
}\sum\limits_{r=0}^{+\infty }\sum\limits_{s=0}^{m}\sum_{l=1}^{d\left(
n,,r,s\right) }\frac{\Phi _{j,p,q}^{m}(z)\overline{\Phi _{l,r,s}^{m}(w)}}{%
\sqrt{\mathcal{N}(z)\mathcal{N}\left( w\right) }}\left\langle \varphi
_{j,p,q},\varphi _{l,r,s}\right\rangle _{\mathcal{H}}  \label{4.15}
\end{equation}
Recalling that
\begin{equation}
\left\langle \varphi _{j,p,q},\varphi _{l,r,s}\right\rangle _{\mathcal{H}%
}=\delta _{j,l}\delta _{p,r}\delta _{q,s}  \label{4.16}
\end{equation}
since $\left\{ \varphi _{j,p,q}\right\} $ is an orthonormal basis of $%
\mathcal{H}$, the above sum in \eqref{4.15} reduces to
\begin{equation}
<z,m|w,m>=\left( \mathcal{N}(z)\mathcal{N}\left( w\right) \right) ^{-\frac{1%
}{2}}\sum\limits_{p=0}^{+\infty }\sum\limits_{q=0}^{m}\sum_{j=1}^{d\left(
n,p,q\right) }\Phi _{j,p,q}^{m}(z)\overline{\Phi _{j,p,q}^{m}(w)}
\label{4.17}
\end{equation}
Now, by Lemma \ref{Lemma3.1}, Eq. \eqref{4.17}  takes the form
\begin{equation}
<z,m|w,m>=\left( \mathcal{N}(z)\mathcal{N}\left( w\right) \right) ^{-\frac{1%
}{2}}\pi ^{-n}e^{\left\langle z,w\right\rangle }L_{m}^{\left( n-1\right)
}\left( \left| z-w\right| ^{2}\right) .  \label{4.18}
\end{equation}
So that the squared modulus of the scalar product in \eqref{4.18}
reads
\begin{equation}
\left| <z,m|w,m>\right| ^{2}=\left( \mathcal{N}(z)\mathcal{N}\left( w\right)
\right) ^{-1}\pi ^{-2n}e^{2{\Re }\left\langle z,w\right\rangle }\left(
L_{m}^{\left( n-1\right) }\left( \left| z-w\right| ^{2}\right) \right) ^{2}.
\label{4.19}
\end{equation}
Returning back to \eqref{4.14} and inserting \eqref{4.19},
we obtain that
\begin{equation}
\mathbb{E}_{\left\{ \mid z,m>\right\} }\left( A_{\varphi }\right) =\int_{%
\mathbb{C}^{n}}\mathcal{N}(z)^{-1}\pi ^{-2n}e^{2{\Re }\left\langle
z,w\right\rangle -|w|^{2}}\left( L_{m}^{\left( n-1\right) }\left( \left|
z-w\right| ^{2}\right) \right) ^{2}\varphi \left( w\right) d\mu \left(
w\right) .  \label{4.20}
\end{equation}
Replacing $\mathcal{N}(z)$ by its expression in Proposition \ref{Proposition4.2.}, we arrive at the result in \eqref{4.10}. This ends the
proof. \hfill $\square $ \newline

Finally, we summarize the above discussion by considering the following
definition.

\begin{definition}
\label{Definition4.2} The Berezin transform of the classical observable $%
\varphi \in L^{2}\left( \mathbb{C}^{n},d\mu \right) $ constructed according
to the quantization by the coherent states $\left\{ \mid z,m>\right\} $ in %
\eqref{4.1} is obtained by associating to $\varphi $ the obtained mean value
in \eqref{4.10}. That is,
\begin{equation}
B_{m}\left[ \varphi \right] (z):=\mathbb{E}_{\left\{ \mid z,m>\right\}
}\left( A_{\varphi }\right)   \label{4.21}
\end{equation}
for every $z\in \mathbb{C}^{n}$.
\end{definition}

\begin{remark}
As mentioned above, for $m=0,$ the eigenspace $A_{0}^{2}(\mathbb{C}^{n})$
coincides with the Bargmann space of analytic functions on $\mathbb{C}^{n}$
that are $e^{-|z|^{2}}d\mu $-square integrable with the reproducing kernel $%
K_{0}(z,w):=\pi ^{-n}e^{\left\langle z,w\right\rangle }$ and the associated
Berezin transform $B_{0}$ is given by a convolution product over the group $%
\mathbb{C}^{n}$ as
\begin{equation}
B_{0}[\phi ](z):=\left( \pi ^{-n}e^{-|w|^{2}}*\phi \right) (z);\qquad \phi
\in L^{2}(\mathbb{C}^{n};d\mu ).  \label{4.22}
\end{equation}
Furthermore, it can be expressed in terms of the Euclidean Laplacian on $%
\mathbb{C}^{n}$ as $B_{0}=e^{\frac{1}{4}\Delta _{\mathbb{C}^{n}}}$ (\cite
{2,20,21}).
\end{remark}

\section{The transform $B_{m}$ and the Euclidean Laplacian}

From \eqref{4.10} and \eqref{4.11} it is easy to see from that the transform

\begin{equation}
B_{m}\left[ \varphi \right] (z)=\frac{m!}{(n)_{m}\pi ^{n}}\int_{\mathbb{C}%
^{n}}e^{-|z-w|^{2}}\left( L_{m}^{\left( n-1\right) }\left( |z-w|^{2}\right)
\right) ^{2}\varphi (w)d\mu (w)  \label{5.1}
\end{equation}
can written as a convolution operator as
\begin{equation}
B_{m}\left[ \varphi \right] =h_{m}\ast \varphi ,\quad \varphi \in
L^{2}\left( \mathbb{C}^{n},d\mu \right) ,  \label{5.2}
\end{equation}
where
\begin{equation}
h_{m}(z)=\frac{m!}{(n)_{m}\pi ^{n}}e^{-|z|^{2}}\left( L_{m}^{\left(
n-1\right) }\left( |z|^{2}\right) \right) ^{2}, \qquad z\in \mathbb{C}^{n}.
\label{5.3}
\end{equation}
In view of \eqref{5.2} a general principle (\cite[p.200]{22}) guaranties
that $B_{m}$ is a function of the Euclidean Laplacian $\Delta _{\mathbb{C}%
^{n}}$. As in \cite{6} we start from the fact that $B_{m}$ should be the
Fourier transform of the function $h_{m}(z)$ evaluated at $\frac{1}{i}$
times the gradient operator $\nabla $. i.e.,
\begin{equation}
B_{m}=\frak{F}\left[ h_{m}\right] (\frac{1}{i}\nabla ).  \label{5.4}
\end{equation}
Here, our method is based on straightforward calculations.

\begin{theorem}
\label{Theorem5.1} Let $m=0,1,2,\cdots $. Then, the Berezin transform $B_{m}$
can be expressed in terms of the Laplacian $\triangle _{\mathbb{C}^{n}}$ as
\begin{equation}
B_{m}=e^{\frac{1}{4}\triangle _{^{_{\mathbb{C}}n}}}\sum\limits_{j=0}^{2m}%
\sigma _{j}^{(n,m)}L_{j}^{(n-1)}\left( \frac{-1}{4}\triangle _{\mathbb{C}%
^{n}}\right)   \label{5.5}
\end{equation}
with
\begin{equation}
\sigma _{j}^{(n,m)}=\frac{m!}{(n)_{m}}\sum\limits_{s=0}^{j}\left(
\begin{array}{l}
j \\
s
\end{array}
\right) \left(
\begin{array}{l}
m+n-1 \\
m-j+s
\end{array}
\right) \left(
\begin{array}{l}
m+n-1 \\
m-s
\end{array}
\right) .  \label{5.6}
\end{equation}
\end{theorem}

\noindent \textit{Proof.} We start by looking at the integral
\begin{equation}
\widehat{h_{m}}(\xi )=\int_{\mathbb{C}^{n}}e^{-i\left\langle \xi
,z\right\rangle }h_{m}(z)d\mu (z).  \label{5.7}
\end{equation}
Inserting \eqref{5.3} in \eqref{5.7} and using polar coordinates $z=\rho
\omega ,\rho >0$ and $\omega \in \mathbb{S}^{2n-1}$, then \eqref{5.7} takes
the form
\begin{align}
\widehat{h}_{m}(\xi )& =\frac{m!}{\pi ^{n}(n)_{m}}\int_{0}^{+\infty
}\int_{S^{2n-1}}e^{-i\left\langle \xi ,z\right\rangle }e^{-\rho ^{2}}\left(
L_{m}^{(n-1)}\left( \rho ^{2}\right) \right) ^{2}\rho ^{2n-1}d\rho d\sigma
(\omega )  \label{5.8} \\
& =\frac{m!}{\pi ^{n}(n)_{m}}\int_{0}^{+\infty }e^{-\rho ^{2}}\left(
L_{m}^{(n-1)}\left( \rho ^{2}\right) \right) ^{2}\rho ^{2n-1}\left(
\int_{S^{2n-1}}e^{-i\left\langle \xi ,\rho \omega \right\rangle }d\sigma
(\omega )\right) d\rho .  \label{5.9}
\end{align}
The last integral in \eqref{5.9} can be identified as a Bochner integral (%
\cite[p.646]{23}), as:
\begin{equation}
\int_{S^{2n-1}}e^{-i2\pi \left\langle \rho \zeta ,\omega \right\rangle
}d\sigma (\omega )=2\pi \rho ^{-n+1}\left| \zeta \right|
^{-n+1}J_{n-1}\left( 2\pi \rho \left| \zeta \right| \right) ,  \label{5.10}
\end{equation}
$J_{\nu }(\cdot )$ being the Bessel function. Therefore, we set $\zeta
=(2\pi )^{-1}\xi $ and we insert \eqref{5.10} into \eqref{5.9} to get that
\begin{equation}
\widehat{h}_{m}(\xi )=\frac{2^{n}m!}{(n)_{m}}|\xi |^{-n+1}\int_{0}^{+\infty
}\rho ^{n}e^{-\rho ^{2}}\left( L_{m}^{(n-1)}\left( \rho ^{2}\right) \right)
^{2}J_{n-1}\left( \rho |\xi |\right) d\rho .  \label{5.11}
\end{equation}
Now, the Feldheim formula \cite{24}, which expresses the product of Laguerre
polynomials as a sum of Laguerre polynomials, is given by
\begin{align}
L_{k}^{(\alpha )}(x)L_{l}^{(\beta )}(x)& =\sum\limits_{j=0}^{k+l}A_{j}\left(
k,l,\alpha ,\beta \right) L_{j}^{(\alpha +\beta )}(x)  \label{5.12} \\
& =\left( -1\right) ^{k+l}\sum\limits_{j=0}^{k+l}A_{j}\left( k,l,\beta
-k+l,\alpha +k-l\right) \frac{x^{j}}{j!}  \label{5.13}
\end{align}
with
\begin{equation}
A_{j}\left( k,l,\alpha ,\beta \right) =\left( -1\right)
^{k+l+j}\sum\limits_{s=0}^{j}\left(
\begin{array}{l}
j \\
s
\end{array}
\right) \left(
\begin{array}{l}
k+\alpha \\
l-j+s
\end{array}
\right) \left(
\begin{array}{l}
l+\beta \\
k-s
\end{array}
\right) ,  \label{5.14}
\end{equation}
$\Re \alpha >-1$, $\Re \beta >-1$, $\Re (\alpha +\beta )>-1$. We make use of
this formula for the particular values of $k=l=m$, $\alpha =\beta =n-1$ and $%
x=\rho ^{2},$ we obtain
\begin{equation}
\left( L_{m}^{(n-1)}(\rho ^{2})\right) ^{2}=\sum\limits_{j=0}^{2m}\gamma
_{j}^{(n,m)}\frac{\rho ^{2j}}{j!}  \label{5.15}
\end{equation}
with
\begin{equation}
\gamma _{j}^{(n,m)}:=(-1)^{j}A_{j}(m,m,n-1,n-1).  \label{5.16}
\end{equation}
Returning back to \eqref{5.11} and replacing \eqref{5.15}, we get
\begin{align}
\widehat{h}_{m}(\xi )& =\frac{2^{n}m!}{(n)_{m}}|\xi
|^{-n+1}\int_{0}^{+\infty }\rho ^{n}e^{-\rho ^{2}}\left(
\sum\limits_{j=0}^{2m}\gamma _{j}^{(n,m)}\frac{\rho ^{2j}}{j!}\right)
J_{n-1}\left( \rho |\xi |\right) d\rho  \label{5.17} \\
& =\frac{2^{n}m!}{(n)_{m}}|\xi |^{-n+1}\sum\limits_{j=0}^{2m}\gamma
_{j}^{(n,m)}\frac{1}{j!}\int_{0}^{+\infty }e^{-\rho ^{2}}\rho
^{2j+n}J_{n-1}(\rho |\xi |)d\rho .  \label{5.18}
\end{align}
Next, we use the identity (\cite[p.704]{19}):
\begin{equation}
\int_{0}^{+\infty }x^{2s+\nu +1}e^{-x^{2}}J_{\nu }\left( 2x\sqrt{z}\right)
dx=\frac{s!}{2}e^{-z}z^{\frac{1}{2}\nu }L_{s}^{\left( \nu \right) }(z);
\label{5.19}
\end{equation}
$s=0,1,2,\cdots $, $s+\Re \mu >-1,$ for $x=\rho $, $\nu =n-1$, $s=j$ and $2%
\sqrt{z}=|\xi |$. This gives
\begin{equation}
\widehat{h}_{m}(\xi )=\frac{m!}{(n)_{m}}e^{-\frac{1}{4}|\xi
|^{2}}\sum\limits_{j=0}^{2m}\gamma _{j}^{(n,m)}L_{j}^{(n-1)}\left( \frac{1}{4%
}|\xi |^{2}\right) .  \label{5.20}
\end{equation}
Finally, we replace $\xi $ by $\frac{1}{i}\nabla $ and we state the first
result as follows. \hfill $\square$\newline

Another way to write the Berezin transform $B_{m}$ as function of Laplacian $%
\triangle _{\mathbb{C}^{n}}$ is as follows.

\begin{theorem}
\label{Theorem5.2} Let $m=0,1,2\cdots $. Then, the Berezin transform $B_{m}$
can be expressed in terms of the Laplacian $\triangle _{\mathbb{C}^{n}}$ as
\begin{equation}
B_{m}=e^{\frac{1}{4}\triangle _{\mathbb{C}^{n}}}\sum\limits_{j=0}^{2m}\kappa
_{j}^{(n,m)}\left( \triangle _{\mathbb{C}^{n}}\right) ^{j}  \label{5.21}
\end{equation}
with
\begin{equation}
\kappa _{j}^{(n,m)}=\frac{2^{2m}(m!)^{3}(-1)^{j}\,{_{3}F_{2}}(\frac{j}{2}-m,%
\frac{j+1}{2}-m,j+n,j-m+1,j-m+1;1)}{(n)_{m}j!2^{3j}(2m-j)!\left( \Gamma
(j-m+1)\right) ^{2}}.  \label{5.22}
\end{equation}
\end{theorem}

\noindent \textit{Proof.} We return back to \eqref{5.11} and we make us of
the following linearization of the product of Laguerre polynomials (%
\cite[p.7361]{25}):
\begin{equation}
L_{k}^{\left( \alpha \right) }(x)L_{l}^{\left( \alpha \right)
}(x)=\sum\limits_{j=\left| k-l\right| }^{k+l}C_{j}\left( k,l,\alpha \right)
L_{j}^{\left( \alpha \right) }(x),  \label{5.23}
\end{equation}
where the coefficients are given in terms of $_{3}\digamma _{2}$
hypergeometric function \cite{23} as
\begin{eqnarray}
C_{j}\left( k,l,\alpha \right)  &=&\frac{2^{k+l-j}k!l!}{\left( k+l-j\right)
!\Gamma \left( j-k+1\right) \Gamma (j-l+1)}  \label{24} \\
&&\times _{3}\digamma _{2}\left( \frac{j-k-l}{2},\frac{j-k-l+1}{2},j+\alpha
+1,j-k+1,j-l+1;1\right)   \notag
\end{eqnarray}
for the particular case $\alpha =n-1,k=l=m$ and $x=\rho ^{2}$. We obtain
\begin{equation}
\left( L_{m}^{(n-1)}\left( \rho ^{2}\right) \right)
^{2}=\sum\limits_{j=0}^{2m}c_{j}^{(n,m)}L_{j}^{(n-1)}\left( \rho ^{2}\right)
,  \label{5.25}
\end{equation}
where
\begin{equation}
c_{j}^{(n,m)}=\frac{2^{2m-j}\left( m!\right) ^{2}._{3}F_{2}\left( \frac{j}{2}%
-m,\frac{j+1}{2}-m,j+n,j-m+1,j-m+1;1\right) }{\left( 2m-j\right) !\left(
\Gamma \left( j-m+1\right) \right) ^{2}}.  \label{5.26}
\end{equation}
Therefore, Eq. \eqref{5.11} takes the form
\begin{align}
\widehat{h}_{m}(\xi )& =\frac{2^{n}m!}{(n)_{m}}|\xi
|^{-n+1}\int_{0}^{+\infty }\rho ^{n}e^{-\rho ^{2}}\left(
\sum\limits_{j=0}^{2m}c_{j}^{(n,m)}L_{j}^{(n-1)}\left( \rho ^{2}\right)
\right) J_{n-1}\left( \rho |\xi |\right) d\rho   \label{5.27} \\
& =\frac{2^{n}m!}{(n)_{m}}|\xi
|^{-n+1}\sum\limits_{j=0}^{2m}c_{j}^{(n,m)}\int_{0}^{+\infty }e^{-\rho
^{2}}\rho ^{n}L_{j}^{(n-1)}\left( \rho ^{2}\right) J_{n-1}\left( \rho |\xi
|\right) d\rho .  \label{5.28}
\end{align}
Next, making use of the identity (\cite[p.812]{19}):
\begin{equation}
\int_{0}^{+\infty }e^{-x^{2}}x^{\nu +1}L_{s}^{\left( \nu \right) }\left(
x^{2}\right) J_{n-1}\left( xu\right) dx=\frac{1}{2s!}\left( \frac{u}{2}%
\right) ^{2s+\nu }e^{-\frac{1}{4}u^{2}}  \label{5.29}
\end{equation}
for $\nu =n-1,x=\rho ,s=j$ and $u=|\xi |,$ the integral in \eqref{5.29}
takes the form
\begin{equation}
\int_{0}^{+\infty }e^{-\rho ^{2}}\rho ^{n}L_{j}^{(n-1)}\left( \rho
^{2}\right) J_{n-1}\left( \rho |\xi |\right) d\rho =\frac{1}{2j!}\left(
\frac{|\xi |}{2}\right) ^{2j+n-1}e^{-\frac{1}{4}|\xi |^{2}}  \label{5.30}
\end{equation}
and \eqref{5.28} becomes
\begin{equation}
\widehat{h}_{m}(\xi )=\frac{m!}{(n)_{m}}e^{-\frac{1}{4}|\xi
|^{2}}\sum\limits_{j=0}^{2m}c_{j}^{(n,m)}\frac{\left( -1\right) ^{j}}{%
j!2^{2j}}\left( -|\xi |^{2}\right) ^{j}.  \label{5.31}
\end{equation}
Finally, we replace $\xi $ by $\frac{1}{i}\nabla $ and we state the
result.\fin

\begin{remark}
\label{Remark5.1} In \cite{6}, we have proved that
\begin{equation}
B_{m}=e^{\frac{1}{4}\triangle _{\mathbb{C}^{n}}}\sum\limits_{k=0}^{m}\frac{%
\left( n-1\right) _{k}\left( m-k\right) !}{(n)_{m}k!}\left( \frac{\Delta _{%
\mathbb{C}^{n}}}{4}\right) ^{k}L_{m-k}^{\left( k\right) }\left( \frac{%
-\Delta _{\mathbb{C}^{n}}}{4}\right) L_{m-k}^{\left( n-1+k\right) }\left(
\frac{-\Delta _{\mathbb{C}^{n}}}{4}\right)   \label{5.32}
\end{equation}
so that formulae in \eqref{5.5} and \eqref{5.21} represent other ways to
write the transform $B_{m}$ as function of $\triangle _{\mathbb{C}^{n}}.$
\end{remark}

\section{Conclusions}

While dealing with a class of generalized Bargmann spaces \cite{3}, we first
have been concerned with a direct proof for obtaining the reproducing kernel
of these spaces starting from the knowledge of an explicit orthonormal basis
and by exploiting an addition formula for Laguerre polynomials involving the
disk polynomials due to Koornwinder \cite{17}. With the help of this basis,
we have constructed for each of these spaces a set of coherent states by
following a generalized formalism in order to apply a coherent state
quantization method \cite{8}. This has provided us with another way to
recover the Berezin transforms attached to the generalized Bargmann spaces
under consideration, which was constructed in \cite{5,6} by means of
Toeplitz operators. For related recent works in the literature, we should
mention the references \cite{26,27,28,29,30}. We have also
established two other formulae expressing these Berezin transforms as
functions of Euclidean Laplacian. These two formulae together with a first
one \cite{6} could be of help in physics problems. Why?. First, we should
note that the Euclidean Laplacian represents (in suitable units) the
Hamiltonian of a \textit{free} particle in quantum mechanics. On the other
hand, in view of Remark 3.1, the transform $B_{m}$ encodes the \textit{effect%
} of the magnetic field at the $m$th eigenenergy (\textit{Landau level}) so
that it could be useful to prepare for physicists formulae expressing $B_{m}$
as function of this Laplacian in all possible different forms. In some
sense, these formulae express, at some energy level, a relation linking a
transform arising from a magnetic Schrodinger operator with a quantity
involving the non-magnetic Schrodinger operator. This link is well defined
through the exponential prefactor (which reflects the free particle case)
times a polynomial function of the Laplacian. The degree and coefficients of
this polynomial function are given explicitly as in Theorem \ref{Theorem5.2}
with the effort to describe the polynomial part in a precise way. The
diamagnetism of spinless bose systems \cite{9} or diamagnetic inequalities
illustrate very well this kind of picture.

\end{document}